\documentclass[12pt]{article}

\usepackage{amsfonts}

\newtheorem{theorem}{Theorem}[section]
\newtheorem{lemma}{Lemma}[section]

\makeatletter
\@addtoreset{equation}{section}
\makeatother

\begin{document}

\begin{center}
{\bf \Large\bf Asymptotics of eigenvalues for an energy operator of the one model of quantum physics.}
\\ \bigskip
E.А.Yanovich
\\ \bigskip
St-Petersburg State Polytechnical University ,\\
  Department of Higher Mathematics,\\
  29 Polytechnicheskaya st., St. Petersburg, 195251 Russia\\
   E-mail: teduard@land.ru
\end{center}

\begin{quote}
{\footnotesize
In this paper we consider eigenvalues asymptotics of the energy operator in the one of the most interesting models of quantum physics, describing an interaction between two-level system and harmonic oscillator. The energy operator of this model can be reduced to some class of infinite Jacobi matrices. Discrete spectrum of this class of operators represents the perturbed spectrum of harmonic oscillator. The perturbation is an unbounded operator compact with respect to unperturbed one. We use slightly modified Janas-Naboko successive diagonalization approach and some new compactness criteria for infinite matrices. Two first terms of eigenvalues asymptotics and the estimation of remainder are found.}
\end{quote}
\bigskip

\section{Introduction and main results.}

We consider the energy operator of the following form
$$
  \hat {\bf H}=\frac{\hbar\omega_0}{2}\,\hat{\bf\sigma}_z+\hbar\omega\:
  \hat {\bf a}^+\,\hat {\bf a}+\hbar\lambda\,
  (\hat{\bf\sigma}_++\hat{\bf\sigma}_-)(\,
  \hat{\bf a}+\hat {\bf a}^+\,)\,,
$$
where $\hat{\bf\sigma}_z,\hat{\bf\sigma}_+,\hat{\bf\sigma}_-$ are the $2\times2$
matrices of form
$$
\hat{\bf\sigma}_z=\left(\begin{array}{cc} 1 & 0\\0 & -1\end{array}\right)
\,,\quad
\hat{\bf\sigma}_+=\left(\begin{array}{cc} 0 & 1\\0 & 0\end{array}\right)
\,,\quad
\hat{\bf\sigma}_-=\left(\begin{array}{cc} 0 & 0\\1 & 0\end{array}\right)\,,
$$
$\hat{\bf a}$ and $\hat{\bf a}^+$ are the creation and annihilation operators for harmonic oscillator, $\lambda$ is the interaction constant, $\omega$ is the oscillator frequency, $\omega_0$ is the transition frequency in the two-level system. These matrices and operators are satisfies by the following commutative relations
$$
[\hat{\bf\sigma}_+,\hat{\bf\sigma}_-]=\hat{\bf\sigma}_z\,,\quad
[\hat{\bf\sigma}_z,\hat{\bf\sigma}_+]=2\,\hat{\bf\sigma}_+\,,\quad
[\hat{\bf\sigma}_z,\hat{\bf\sigma}_-]=-2\,\hat{\bf\sigma}_-\,,\quad
[\hat{\bf a},\hat{\bf a}^+]=1
$$

In the work~\cite{1} it was shown that the hamiltonian of this model is represented by the two Jacobi matrices. These matrices have the following general form
\begin{equation}
\label{model}
A=\left(\begin{array}{cccccc} c_1&g\sqrt{1}&0&0&0&\ldots\\
                      g\sqrt{1}&1+c_2&g\sqrt{2}&0&0&\ldots\\
                      0&g\sqrt{2}&2+c_1&g\sqrt{3}&0&\ldots\\
                      0&0&g\sqrt{3}&3+c_2&g\sqrt{4}&\ldots\\
                      0&0&0&g\sqrt{4}&4+c_1&\ldots\\
                      \ldots&\ldots&\ldots&\ldots&\ldots&\ldots
  \end{array}\right)
\end{equation}
where $g$, $c_1$, $c_2$ are real parameters. As it is well known~\cite{2,3} that the matrix $A$ defines a selfadjoint operator with simple spectrum and the domain $D(A)$ is dense in the space $l_2(\mathbb N)$. As the operator $A$ can be considered as relatively compact perturbation of the main diagonal, its spectrum is discrete.

The main goal of this paper is the investigation of the eigenvalues $\lambda_n(A)$ behavior for large values of $n$ and fixed other parameters. The result is given by the following asymptotic formula (Theorem~(\ref{Jaynes}))
$$
\lambda_n(A)=n-g^2+\frac{c_1+c_2}{2}+O\left(\frac{1}{n^{1/16}}\right)\,,
\quad n\to\infty\quad(g\ne 0)
$$

\section{Selection of the main component in the asymptotics.}

Let's present the operator $A$ in~(\ref{model}) in the form
\begin{equation}
\label{r}
A=A_0+\frac{c_1+c_2}{2}\,I+\frac{c_1-c_2}{2}\,R\,,
\end{equation}
where $I$ is the identical matrix, $A_0$ and $R$ are defined in the following way
$$
A_0=\left(\begin{array}{ccccc} 0&g\sqrt{1}&0&0&\ldots\\
                      g\sqrt{1}&1&g\sqrt{2}&0&\ldots\\
                      0&g\sqrt{2}&2&g\sqrt{3}&\ldots\\
                      0&0&g\sqrt{3}&3&\ldots\\
                      \ldots&\ldots&\ldots&\ldots&\ldots
  \end{array}\right)\,,\quad
R=\left(\begin{array}{cccccc} 1&0&0&0&\ldots\\
                      0&-1&0&0&\ldots\\
                      0&0&1&0&\ldots\\
                      0&0&0&-1&\ldots\\
                      \ldots&\ldots&\ldots&\ldots&\ldots
  \end{array}\right)
$$

The matrix $A_0$ represents so called "shifted oscillator" operator:
$a^+a+g\,(a+a^+)$, where $a^+$ и $a$ are the creation and annihilation operators. If we use its matrix representation, we obtain exactly the matrix $A_0$.

Eigenvalues problem for the operator $A_0$ has an exact solution. This solution can be obtained by the different ways. For example, with the help of Bogolubov's transformation \cite{4} or by using continued fractions~\cite{5}. This solution has the form
\begin{equation}
\label{Bogolubov1}
A_0\,e_n=\mu_n\,e_n\,,\quad \mu_n=n-g^2\,,\quad n=0,1,2,\ldots
\end{equation}
Here $e_n$ are the normalized eigenvectors of the operator $A_0$. Its expansion through
the basis vectors $a_n$ of the matrix representation~(\ref{model}) has the form
\begin{equation}
\label{Bogolubov2}
e_m=Ua_m=\sum_{n=0}^\infty U_{n,m}\,a_n\,,\quad U_{n,m}=
\exp\{-g^2/2\}\,\sqrt{\frac{n!}{m!}}\:\,g^{m-n}\,L_n^{(m-n)}(g^2)=
\omega_n^{(m-n)}(g^2)\,,
\end{equation}
where the matrix elements $U_{n,m}$ of the orthogonal transformation $U$ are defined
by Feynman-Schwinger's formula~\cite{9,10} (see also~\cite{11});
$L_n^{(s)}(x)$ are the generalized Chebyshev-Laguerre's polynomials~\cite{13}
$$
L_n^{(s)}(x)=\frac{(n+s)!}{n!}\,\sum_{i=0}^n C_n^i\,(-1)^i\,\frac{x^i}{(i+s)!}
\,,\quad C_n^i=\frac{n!}{i!\,(n-i)!}\,,\quad (s\ge 0)
$$
$$
L_n^{(-s)}(x)=(-x)^s\,\frac{(n-s)!}{n!}\,L_{n-s}^{(s)}(x)\,,
$$
and $\omega_n^{(s)}(x)$ are the normalized Laguerre's functions
$$
\omega_n^{(s)}(x)=\sqrt{\frac{n!}{(n+s)!}}\:\,e^{-x/2}\,x^{s/2}\,L_n^{(s)}(x)
\,,\quad
\int\limits_0^{+\infty}\omega_m^{(s)}(x)\,\omega_n^{(s)}(x)\,dx=\delta_{n,m}
$$

Let's note here that the comleteness of the "shifted oscillator" eigenfunctions (in coordinate representation) for complex values of the parameter $g$ is considered in ~\cite{12}.

Let's us find the matrix of the operator $R$ in the basis of the operator $A_0$ eigenvectors. Denoting the elements of the transformed matrix as $\tilde R_{k,m}$ and taking into account that $R_{n,m}=(-1)^n\,\delta_{n,m}$, we obtain
\begin{equation}
\label{R-transformed}
\tilde R_{k,m}=(Re_m,e_k)=(U^TRU)_{k,m}=\sum_{n=0}^\infty (-1)^n\,U_{n,k}\,
U_{n,m}
\end{equation}
Let's present the matrix elements $U_{n,m}$ as contour integral
\begin{equation}
\label{contour}
U_{n,m}=
\exp\{-g^2/2\}\,\sqrt{\frac{m!}{n!}}\:\,g^{n-m}\,\frac{1}{2\pi i}\,
\oint\limits_C x^{m-1}\left(\frac{1}{x}-1\right)^n
\exp\left\{\frac{g^2}{x}\right\}dx\,,
\end{equation}
where $C$ is a unit circle centered in the origin of the complex plane $x$.

Substituting~(\ref{contour}) in~(\ref{R-transformed}) and summing up over $n$, we find
$$
\tilde R_{k,m}=
\exp\{-2g^2\}\,\sqrt{k!\,m!}\;\,g^{-m-k}\,
\frac{1}{(2\pi i)^2}\,\oint\limits_C\oint\limits_C (x)^{m-1}\,(x')^{k-1}
\,\exp\left\{g^2\left(\frac{2}{x}+\frac{2}{x'}-\frac{1}{xx'}
\right)\right\}dx\,dx'
$$

Contour integrals in this expression can be calculated consistently with the help of residues. As a result, we obtain
$$
\tilde R_{k,m}=(-1)^k\,\exp\{-2g^2\}\,\sqrt{\frac{m!}{k!}}\:\,(2g)^{m-k}\,
\sum_{i=0}^k C_k^i\,(-1)^i\,\frac{(4g^2)^i}{(i+m-k)!}
$$
(if the expression $(i+m-k)$ in the denominator of the last sum is negative, the corresponding term must be considered equal to zero).

Comparing this expression with~(\ref{Bogolubov2}), we obtain
\begin{equation}
\label{R-matrix}
\tilde R_{k,m}=(-1)^k\,\omega_k^{(m-k)}(4g^2)
\end{equation}

In spite of seeming asymmetry, the matrix $\tilde R_{k,m}$ is symmetric ($\tilde R_{k,m}=\tilde R_{m,k}$). It can be easily verified with the help of the known property of the generalized Chebyshev-Laguerre's polynomials~\cite{13}
$$
L_n^{(-s)}(x)=(-x)^s\,\frac{(n-s)!}{n!}\,L_{n-s}^{(s)}(x)
$$

Using the asymtotics of the generalized Chebyshev-Laguerre's polynomials~\cite{13}
\begin{equation}
\label{Laguerre_asym}
L_{n}^{s}(x)=\pi^{-1/2}\,n^{s/2-1/4}\,x^{-s/2-1/4}\,e^{x/2}\left\{
\cos(2\sqrt{nx}-s\pi/2-\pi/4)+O(n^{-1/2})\right\}\,,\:n\to\infty\,,
\end{equation}
we find
\begin{equation}
\label{Condition_in_Theorem}
\lim_{n\to\infty} \tilde R_{n,n+p} =0\,,\quad\forall p\in Z
\end{equation}

In what follows we will need the following result~\cite{6}

\begin{lemma}[J.Janas-S.Naboko]
\label{Janas-Naboko}
Let $D$ be a selfadjoint operator in a Hilbert space $H$ with simple discrete spectrum ($De_n=\mu_n e_n$), where $\{e_n\}$ is an orthonormal basis of eigenvectors in $H$ and $\mu_n$ are simple eigenvalues ($\mu_n\to\infty$), ordered by $|\mu_i|\le|\mu_{i+1}|$.  Assume that $|\mu_i-\mu_k|\ge\epsilon_0>0\,,\,\,\forall i\ne k$. If $R$ is a compact operator in $H$ then the eigenvalues $\lambda_n(T)$ of the operator $T=D+R$ (with discrete spectrum too) become simple for large values of $n$ and satisfy to the asymptotic formula
\begin{equation}
\label{Naboko}
\lambda_n(T)=\mu_n+O(\|R^*e_n\|)\,,\quad n\to\infty\,,
\end{equation}
where $R^*$ is the adjoint operator with respect to $R$.
\end{lemma}

The matrix $R$ and the matrix similar to it $\tilde R=U^TRU$, represents bounded noncompact operator (projector) since $R^2=I$. Therefore we can't apply here at once the Lemma~(\ref{Janas-Naboko}).

Let's prove the following theorem:

\begin{theorem}
\label{theorem2}
Let $D$ be a selfadjoint operator in a Hilbert space $H$ with eigenvalues $\mu_n=n$, ($n=0,1,2,\ldots$) and complete system of corresponding eigenvectors in $H$. Let $R$ be a bounded, selfadjoint, noncompact operator and its matrix $R_{n,k}$ in the basis of operator $D$ eigenvectors satisfy to the condition
\begin{equation}
\label{t2_1}
\lim_{n\to\infty} R_{n,\,n+p}=0\,,\quad\forall p\in Z
\end{equation}
Then the eigenvalues $\lambda_n(T)$ of the operator $T=D+R$ (having a discrete spectrum too) become simple for large values of $n$ and satisfy to the following asymptotic estimation
\begin{equation}
\label{t2_2}
\lambda_n(T)=n+R_{n,n}+O\left(s_n\right)\,,
\quad n\to\infty\,,
\end{equation}
where
$$
s_n=\sqrt{\sum_{k\ne n}\frac{|R_{k,n}|^2}{(n-k)^2}}\,,
$$
and $s_n\to0$ at $n\to\infty$.
\end{theorem}

For the proof of this theorem we need the following compactness criteria for infinite matrices.

\begin{lemma}
\label{compactness}
Let $V$ be a bounded, noncompact operator in a Hilbert space $H$. Let its matrix $V_{i,j}$ ($i,j=0,1,\ldots$) in some orthonormal basis satisfy to the condition
\begin{equation}
\label{condition_1}
\lim_{n\to\infty} V_{n,\,n+p}=0\,,\quad\forall p\in Z
\end{equation}
Let $b=\{b_i\}_{i=-\infty}^\infty$ is an arbitrary $l_2$-sequence
\begin{equation}
\label{condition_2}
\|b\|^2=\sum_{i=-\infty}^\infty |b_i|^2<\infty
\end{equation}
Then the operator $K$ with matrix $K_{i,j}=b_{i-j}V_{i,j}$ ($i,j=0,1,\ldots$) is compact in $H$.
\end{lemma}

{\bf Proof.}\, Let's show at first that the operator $K$ is bounded. For that we need to prove the estimates~\cite{7}
\begin{equation}
\label{bounded}
\sum_{j=0}^\infty|K_{i,j}|<A\,,\:\:\forall i\,;\quad\quad
\sum_{i=0}^\infty|K_{i,j}|<A\,,\:\:\forall j\,,
\end{equation}
where $A$ is a constant independent of $i$ and $j$ ($\|K\|\le A$). Using Cauchy's inequality, we have
$$
\sum_{j=0}^\infty|K_{i,j}|=\sum_{j=0}^\infty|b_{i-j}V_{i,j}|\le
\left(\sum_{j=0}^\infty |b_{i-j}|^2\right)^{1/2}
\left(\sum_{j=0}^\infty |V_{i,j}|^2\right)^{1/2}\le
\|b\|\sqrt{(VV^*)_{i,i}}\le \|b\|\cdot\|V\|
$$
Due to~(\ref{condition_2}) the first estimate in~(\ref{bounded}) is fulfilled. By the same way the validity of the second estimatein~(\ref{bounded}) is established. Thus the operator $K$ is compact. Let's prove now its compactness.

Let's define the cut-off function $b^{(n)}=\{b^{(n)}_i\}_{i=-\infty}^\infty$ of the sequence $\{b_i\}$
$$
b_i^{(n)}=\left\{
                  \begin{array}{lc}
                  0\,,&|i|>n\\
                  b_i\,,&|i|\le n\\
                  \end{array}
          \right.
$$

Let's define the sequence of operators $K^{(n)}$ by the formula $K^{(n)}_{i,j}=b^{(n)}_{i-j}V_{i,j}$. It follows from this definition and from ~(\ref{condition_1}) that $K^{(n)}$ is a compact operator for arbitrary $n$. We have further so as at the proof of~(\ref{bounded})
$$
\|K-K^{(n)}\|\le \|b-b^{(n)}\|\cdot\|V\|
$$
Terefore
$$
\|K-K^{(n)}\|\to 0\,,\quad n\to\infty\,,
$$
and $K$ is compact as a limit by norm of compact operators~\cite{7,8}.\\
-------------------------------------------------

{\bf Proof of the theorem~(\ref{theorem2}).}\,
Let's associate to each operator the matrix in the basis of the operator $D$ eigevectors, remaining the same notations.

Following the main ideas of the work~\cite{6} let's show that there exist such anti-hermitian operator $K$ ($K^*=-K$) that
\begin{equation}
\label{equivalent}
(I+K)T-D_1(I+K)=B\,,
\end{equation}
where $B$ is compact operator and $D_1=D+\mbox{diag}\{R_{n,n}\}$.
(So $D_1$ is the diagonal matrix with elements $(D_1)_{n,n}=n+R_{n,n}$). Suppose that such operator $K$ have found. The condition~(\ref{equivalent}) means that
$$
T=(I+K)^{-1}(D_1+B(I+K)^{-1})(I+K)
$$
(The existence of the inverse operator $(I+K)^{-1}$ follows from the anti-hermitianess of $K$). That is the operators $T$ and $D_1+B(I+K)^{-1}$ are similar and have the same spectrum. But the operator $B$ is compact and the eigenvalues of $D_1$ due to~(\ref{t2_1}) satisfy to the requirements of Lemma~(\ref{Janas-Naboko}). Applying  the Lemma~(\ref{Janas-Naboko}) we obtain
\begin{equation}
\label{t2_3}
\lambda_n(T)=n+R_{n,n}+O\left(\|B^*e_n\|\right)
\end{equation}
Therefore for the proof of the theorem we should establish the existence of such the operator $K$ and find the matrix of the compact operator $B$. Substituting the expressions for the matrices $T$ и $D_1$ in~(\ref{equivalent}) we obtain
\begin{equation}
\label{t2_4}
(I+K)T-D_1(I+K)=R_1-[D,K]+KR-\mbox{diag}\{R_{n,n}\}\,K\,,
\end{equation}
where $[\cdot,\cdot]$ is the commutator and $R_1$ is the matrix of the operator $R$ with zero main diagonal ($R_1=R-\mbox{diag}\{R_{n,n}\}$).

This expression will be the matrix of compact operator if we can find such an compact operator $K$ that the condition $[D,K]=R_1$ is valid, or in matrix form: $K_{i,j}\,(i-j)=(R_1)_{i,j}$. It follows from that
\begin{equation}
\label{t2_5}
K_{i,j}=\frac{R_{i,j}}{i-j}\,,\quad i\ne j\,;\quad\quad K_{i,i}=0\,,\:\:
i=0,1,\ldots
\end{equation}

As the operator $R$ is selfadjpint the corresponding to the matrix~(\ref{t2_5}) operator $K$ is anti-hermitian. Its compactness follows from the Lemma~(\ref{compactness}). Actually, if we choose the sequence $\{b_i\}$ as $\{1/i\}$ ($i\ne0$) then from~(\ref{t2_1}) it follows that all conditions of the Lemma~(\ref{compactness}) are fulfilled.

Now from~(\ref{t2_4}) we find the form of the compact operator $B$:
$$
B=KR-\mbox{diag}\{R_{n,n}\}\,K\,,
$$
and taking into account~(\ref{t2_5}) we obtain
$$
O\left(\|B^*e_n\|\right)=O\left(\|K^*e_n\|\right)=O\left(\sqrt{\sum_{k\ne n}
\frac{|R_{k,n}|^2}{(n-k)^2}}\right)
$$

Substituting this estimate to~(\ref{t2_3}), we obtain the formula~(\ref{t2_2}). The theorem is proved.\\
-------------------------------------------------

Due to~(\ref{Condition_in_Theorem}), the condition~(\ref{t2_1}) of the Theorem~(\ref{theorem2}) is fulfilled. Hence, applying the Theorem~(\ref{theorem2}) and taking into account~(\ref{Bogolubov1}) and~(\ref{r}) we have the following result
$$
\lambda_n(A)=n-g^2+\frac{c_1+c_2}{2}+O\left(s_n\right)\,,\quad n\to\infty\,,
$$
where
\begin{equation}
\label{rem_term}
s_n=\sqrt{\sum_{k\ne n}\frac{|\tilde R_{k,n}|^2}{(n-k)^2}}=
\sqrt{\sum_{k\ne n}\frac{|\omega_k^{(n-k)}(4g^2)|^2}{(n-k)^2}}
\end{equation}

\section{Estimation of the remainder.}

To estimate the decreasing rate of the sequence $s_n$ we should have another estimation for the Laguerre's functions $\omega_n^{(s)}(x)$ than the estimation following from~(\ref{Laguerre_asym}) (in~(\ref{Laguerre_asym}) the parameter $s$ is fixed). We could not find this result among known one and therefore we give here not only the formulation but the proof of it.

\begin{lemma}
\label{unequality1}
Suppose that $x>0$, $s\in Z_+$. Then the following estimate for the Bessel functions $J_s(x)$ is valid
$$
\left|J_s(x)\right|\le 2\,\sqrt{\frac{2}{\pi x}}\,\left(1+
\frac{s}{x}\right)^s
$$
\end{lemma}

{\bf Proof.}\, Let's use known representation~\cite{18}
$$
J_s(x)=\sqrt{\frac{2}{\pi x}}\,\left(P(x,s)\,\cos(x-s\pi/2-\pi/4)-Q(x,s)\,
\sin(x-s\pi/2-\pi/4)\right)\,,
$$
where
$$
P(x,s)=\frac{1}{2\,\Gamma(s+1/2)}\,\int\limits_0^\infty e^{-u}\,u^{s-1/2}
\left\{\left(1+\frac{iu}{2x}\right)^{s-1/2}+\left(1-\frac{iu}{2x}\right)
^{s-1/2}\right\}\,du
$$
$$
Q(x,s)=\frac{1}{2i\,\Gamma(s+1/2)}\,\int\limits_0^\infty e^{-u}\,u^{s-1/2}
\left\{\left(1+\frac{iu}{2x}\right)^{s-1/2}-\left(1-\frac{iu}{2x}\right)
^{s-1/2}\right\}\,du
$$
It is evident that
\begin{equation}
\label{l2_1.1}
|J_s(x)|\le\sqrt{\frac{2}{\pi x}}\,\left(|P(x,s)|+|Q(x,s)|\right)\,,
\end{equation}
and everything reduces to the estimation of the integrals $P(x,s)$ and $Q(x,s)$. Let us consider the integral for $P(x,s)$. The estimation for $Q(x,s)$ is the same.
At $s=0$ we have $|P(x,s)|\le 1$, $|Q(x,s)|\le 1$ and
the estimation~(\ref{l2_1.1}) gives the required inequality. Suppose that $s\in N$.
In this case we have
$$
|P(x,s)|\le\frac{1}{\Gamma(s+1/2)}\,\int\limits_0^\infty e^{-u}\,u^{s-1/2}\,
\left(1+\frac{u}{2x}\right)^{s-1/2}\,du
$$
Expanding the binomial in this integral in the series on $u/2x$
$$
\left(1+\frac{u}{2x}\right)^{s-1/2}=1+\sum_{k=1}^{p-1}
\frac{(s-1/2)\cdot\ldots\cdot(s-1/2-(k-1))}{k!}\,\left(\frac{u}{2x}\right)^k+
$$
$$
+\frac{(s-1/2)\cdot\ldots\cdot(s-1/2-(p-1))}{k!}\,(1+\theta)^{s-p-1/2}\,
\left(\frac{u}{2x}\right)^p\,,\quad \theta\in(0,u/2x)
$$
and putting $p=s$ we have $(1+\theta)^{s-p-1/2}<1$ and therefore
$$
\left(1+\frac{u}{2x}\right)^{s-1/2}<1+\sum_{k=1}^{s}
\frac{(s-1/2)\cdot\ldots\cdot(s-1/2-(k-1))}{k!}\,\left(\frac{u}{2x}\right)^k
$$
Integrating by terms we obtain
$$
|P(x,s)|\le 1+\sum_{k=1}^{s}\frac{\Gamma(s+k+1/2)}{\Gamma(s+1/2)}\,\,
\frac{(s-1/2)\cdot\ldots\cdot(s-1/2-(k-1))}{k!}\,\frac{1}{(2x)^k}<
$$
$$
<1+\sum_{k=1}^{s}(2s)^k\,
\frac{s\cdot\ldots\cdot(s-(k-1))}{k!}\,\frac{1}{(2x)^k}=
\left(1+\frac{s}{x}\right)^s
$$
For $Q(x,s)$ the same estimate is valid and the formula~(\ref{l2_1.1}) leads again to the required inequality. Lemma is proved.\\
-------------------------------------------------

\begin{lemma}
\label{unequality2}
If $x>0$; $n,s\in Z_+$ and $s^{16}\le n$ then
\begin{equation}
\label{main_unequality}
\left|\omega_n^{(s)}(x)\right|\le\frac{C(x)}{(n+1)^{1/4}}\,,
\end{equation}
where the constant $C(x)$ depends on $x$ only.
\end{lemma}

{\bf Proof.}\, Let's use Laguerre's functions integral representation trough the Bessel functions~\cite{13}
$$
\omega_n^{(s)}(x)=\frac{\textstyle e^{x/2}}{\sqrt{n!\,(n+s)!}}
\int\limits_0^\infty e^{-t}\,t^{n+\frac{\scriptstyle s}{2}}\,
J_s(2\sqrt{tx})\,dt\,,\quad n,s\in Z_+
$$
Let us split this integral into two one
$$
\omega_n^{(s)}(x)=\frac{\textstyle e^{x/2}}{\sqrt{n!\,(n+s)!}}
\left(\int\limits_0^{t_0}\:+\:\int\limits_{t_0}^\infty\right)\,,
$$
where $t_0\ge 0$ is an arbitrary now.

For estimation of the first integral let's use the known inequality~\cite{18}
$$
\left|J_s(x)\right|\le 1\,,\quad s\in Z_+\,,\:\:x\in R\,,
$$
For estimation of the second integral we use more precise estimate from Lemma~(\ref{unequality1})
$$
\left|J_s(x)\right|\le 2\,\sqrt{\frac{2}{\pi x}}\,\left(1+
\frac{s}{x}\right)^s<2e\,\sqrt{\frac{2}{\pi x}}\,,\quad x\ge s^2
$$
Putting $t_0=s^4/4x$ (so that at $t\ge t_0$ one can use the last estimate) we have
$$
|\omega_n^{(s)}(x)|\le\frac{\textstyle e^{x/2}}{\sqrt{n!\,(n+s)!}}
\left[\int\limits_0^{t_0} e^{-t}\,t^{n+\frac{\scriptstyle s}{2}}\,dt+
\frac{2e}{\sqrt{\pi\sqrt x}}\,
\int\limits_{t_0}^\infty e^{-t}\,t^{n+\frac{\scriptstyle s}{2}-1/4}\,dt\right]
\le
$$
$$
\le\frac{\textstyle e^{x/2}}{\sqrt{n!\,(n+s)!}}
\left[t_0\:\max_{t\ge0}\left\{e^{-t}\,t^{n+\frac{\scriptstyle s}{2}}\right\}+
\frac{2e}{\sqrt{\pi\sqrt x}}\,\Gamma(n+s/2+3/4)\right]=
$$
$$
=e^{x/2}\,
\left[\frac{s^4}{4x}\:\frac{e^{-(n+\frac{\scriptstyle s}{2})}\,
(n+s/2)^{n+\frac{\scriptstyle s}{2}}}{\sqrt{n!\,(n+s)!}}+
\frac{2e}{\sqrt{\pi\sqrt x}}\,\frac{\Gamma(n+s/2+3/4)}{\sqrt{n!\,(n+s)!}}
\right]
$$

Using known inequalities for Gamma-function following from Stirling's formula
$$
C_1\,z^{z-1/2}\,e^{-z}\le\Gamma(z)\le C_2\,z^{z-1/2}\,e^{-z}\,,\quad
z\ge\delta>0\,,
$$
where $C_1,C_2$ are some constants independent of $z$, let's estimate each term in square brackets. We have
$$
\frac{2e}{\sqrt{\pi\sqrt x}}\,\frac{\Gamma(n+s/2+3/4)}{\sqrt{n!\,(n+s)!}}\le
C(x)\,\frac{(n+s/2+3/4)^{n+s/2+1/4}}{(n+1)^{n/2+1/4}\:(n+s+1)^{n/2+s/2+1/4}}
\le
$$
$$
\le\frac{C(x)}{(n+1)^{1/4}}\,\frac{(1+\frac{s}{2(n+1)})^n}
{(1+\frac{s}{n+1})^{n/2}}\le\frac{C(x)}{(n+1)^{1/4}}\,\left(1+\frac{s^2}
{4(n+1)^2}\right)^{n/2}\le\frac{C(x)}{(n+1)^{1/4}}\,
\exp\left\{\frac{s^2}{8(n+1)}\right\}
$$

At last, if we put $n\ge s^2$ then
\begin{equation}
\label{l2.2.1}
\frac{2e}{\sqrt{\pi\sqrt x}}\,\frac{\Gamma(n+s/2+3/4)}{\sqrt{n!\,(n+s)!}}\le
\frac{C(x)}{(n+1)^{1/4}}\,,\quad (n\ge s^2)
\end{equation}
Similarly, we can estimate the second term
$$
\frac{s^4}{4x}\:\frac{e^{-(n+\frac{\scriptstyle s}{2})}\,
(n+s/2)^{n+\frac{\scriptstyle s}{2}}}{\sqrt{n!\,(n+s)!}}\le
C(x)\,\frac{s^4\,(n+s/2)^{n+s/2}}{(n+1)^{n/2+1/4}\:(n+s+1)^{n/2+s/2+1/4}}\le
$$
$$
\le\frac{C(x)}{(n+1)^{1/4}}\,\frac{s^4}{(n+s+1)^{1/4}}\,
\frac{(1+\frac{s}{2(n+1)})^n}{(1+\frac{s}{n+1})^{n/2}}\le
\frac{C(x)}{(n+1)^{1/4}}\,\frac{s^4}{(n+s+1)^{1/4}}\,
\exp\left\{\frac{s^2}{8(n+1)}\right\}
$$
If $n\ge s^{16}$ then
\begin{equation}
\label{l2.2.2}
\frac{s^4}{4x}\:\frac{e^{-(n+\frac{\scriptstyle s}{2})}\,
(n+s/2)^{n+\frac{\scriptstyle s}{2}}}{\sqrt{n!\,(n+s)!}}\le
\frac{C(x)}{(n+1)^{1/4}}\,,\quad (n\ge s^{16})
\end{equation}
From~(\ref{l2.2.1}),~(\ref{l2.2.2}) it follows that
$$
\left|\omega_n^{(s)}(x)\right|\le\frac{C(x)}{(n+1)^{1/4}}\,,\quad
(n\ge s^{16})\,,
$$
q.e.d.\\
-------------------------------------------------

Having the estimate~(\ref{main_unequality}) and the orthogonality condition of the transformation $U$:
\begin{equation}
\label{cons_of_orthog}
\sum_{k=0}^\infty|U_{k,n}|^2=\sum_{k=0}^\infty\left|\omega_k^{(n-k)}\right|^2=
\sum_{k=0}^n\left|\omega_{n-k}^{(k)}\right|^2+
\sum_{k=1}^\infty\left|\omega_n^{(k)}\right|^2=1\,,\quad\forall n\in Z_+
\end{equation}
( $|\omega_{k}^{(n-k)}|=|\omega_{n}^{(k-n)}|$ ) one can give the estimate of the remainder which is defined by the sum
\begin{equation}
\label{sum4}
\sum_{k\ne n}\frac{|\omega_k^{(n-k)}|^2}{(n-k)^2}=
\sum_{k=1}^n\frac{|\omega_{n-k}^{(k)}|^2}{k^2}+
\sum_{k=1}^\infty\frac{|\omega_n^{(k)}|^2}{k^2}
\end{equation}
From Lemma~(\ref{unequality2}) it follows that
$$
|\omega_{n-k}^{(k)}|^2\le\frac{C}{(n-k+1)^{1/2}}\,,\quad n-k\ge k^{16}\:\:\:
(n\ge k^{16}+k)
$$
$$
|\omega_{n}^{(k)}|^2\le\frac{C}{(n+1)^{1/2}}\,,\quad n\ge k^{16}\:\:\:
(k\le n^{1/16})
$$

Let $k_n\ge 0$ be a maximal nonnegative integer of $k$, satisfying to the equation $n\ge k^{16}+k$. It is evident that $k_n\le n^{1/16}$. Hence
$$
|\omega_{n-k}^{(k)}|^2\le\frac{C}{(n-n^{1/8}+1)^{1/2}}\,,\quad k\le k_n
$$
$$
|\omega_{n}^{(k)}|^2\le\frac{C}{(n+1)^{1/2}}\,,\quad k\le k_n
$$
Let's present the sum~(\ref{sum4}) in the form
$$
\sum_{k\ne n}\frac{|\omega_k^{(n-k)}|^2}{(n-k)^2}=
\left[\sum_{k=1}^{k_n}\frac{|\omega_{n-k}^{(k)}|^2}{k^2}+
\sum_{k=1}^{k_n}\frac{|\omega_n^{(k)}|^2}{k^2}\right]+
\left[\sum_{k=k_n+1}^{n}\frac{|\omega_{n-k}^{(k)}|^2}{k^2}+
\sum_{k=k_n+1}^\infty\frac{|\omega_n^{(k)}|^2}{k^2}\right]
$$
Due to last inequalities we have
$$
\left[\sum_{k=1}^{k_n}\frac{|\omega_{n-k}^{(k)}|^2}{k^2}+
\sum_{k=1}^{k_n}\frac{|\omega_n^{(k)}|^2}{k^2}\right]\le
\frac{2C}{(n-n^{1/8}+1)^{1/2}}\,\sum_{k=1}^{k_n}\frac{1}{k^2}=
O\left(\frac{1}{n^{1/2}}\right)
$$
Since $k_n\sim n^{1/16}$, we have using~(\ref{cons_of_orthog})
$$
\left[\sum_{k=k_n+1}^{n}\frac{|\omega_{n-k}^{(k)}|^2}{k^2}+
\sum_{k=k_n+1}^\infty\frac{|\omega_n^{(k)}|^2}{k^2}\right]\le
\frac{1}{(k_n+1)^2}\,\left[\sum_{k=k_n+1}^{n}|\omega_{n-k}^{(k)}|^2+
\sum_{k=k_n+1}^\infty|\omega_n^{(k)}|^2\right]\le
$$
$$
\le\frac{1}{(k_n+1)^2}=
O\left(\frac{1}{n^{1/8}}\right)
$$
Combining both estimates, we obtain
$$
\sum_{k\ne n}\frac{|\omega_k^{(n-k)}|^2}{(n-k)^2}=
O\left(\frac{1}{n^{1/2}}\right)+O\left(\frac{1}{n^{1/8}}\right)=
O\left(\frac{1}{n^{1/8}}\right)
$$
Taking into account the formula~(\ref{rem_term}), we come to the following main result
\begin{theorem}
\label{Jaynes}
The eigenvalues $\lambda_n(A)$ of the operator $A$~(\ref{model}) at $g\ne 0$ have the following asymptotics
$$
\lambda_n(A)=n-g^2+\frac{c_1+c_2}{2}+O\left(\frac{1}{n^{1/16}}\right)\,,
\quad n\to\infty
$$
\end{theorem}
-------------------------------------------------

{\bf Acknowledgment.} I am grateful to Prof. S.N.Naboko for the discussion of the results. Special gratitude I express to Prof. B.S.Mityagin for his useful remarks.


\begin{thebibliography}{14}
\bibitem{1}
Tur E.A., \textit{Jaynes-Cummings model: Solution without rotating wave
approximation} , Optics and Spectroscopy, Vol. 89, n. 4, 2000, pp. 574-588.
\bibitem{2}
Akhiezer N., Krein M., \textit{Some questions in the theory of moments}, Translation of Mathematical Monographs, Vol. 2, Amer. Math. Soc., Providence, RI, 1962.
\bibitem{3}
Akhiejzer N., \textit{The classical moment problem}, New-York:
Hafner, 1965.
\bibitem{4}
Nagy K., \textit{State Vector Spaces with Indefinite Metric in Quantum Field Theory}, P. Noordhoff, Groningen, Netherlands, 1966.
\bibitem{5}
Tur E.A., \textit{Solution of quntum mecanical problems with th ehelp of continued fractions}, PHD thesis, St-Petersburg State University, 2002.
\bibitem{9}
Feynman R.P., Phys.Rev., {\bf 84}, 1951, 108.
\bibitem{10}
Schwinger J., Phys.Rev., {\bf 91}, 1953, 728.
\bibitem{11}
Baz A., Zel'dovich Y., Perelomov A., \textit{Scattering, reactions and decays in non- relativistic quantum mechanics}, Moscow, 1971.
\bibitem{12}
Kato T., \textit{Perturbation theory for linear operators}, Springer-Verlag Berlin
$\cdot$Heidelberg$\cdot$New York, 1966.
\bibitem{13}
Szego G., \textit{Orthogonal polynomials}, New York, 1939.
\bibitem{6}
Janas J., Naboko S., \textit{Infinite Jacobi matrices with unbounded entries:
Asymptotics of eigenvalues and the transformation operator approach},
SIAM Journal on Mathematical Analysis, 2004, v. 36, n. 2, pp. 643-658.
\bibitem{7}
Smirnov V., \textit{A course of higher mathematics}, v. 5, Moskow, 1960.
\bibitem{8}
Akhiejzer N., Glazman I., \textit{The theory of linear operators
in Hilbert space}, Science, Moscow, 1966.
\bibitem{18}
Watson G., \textit{The theory of Bessel functions}, Cambridge, University Press, 1922.
\end{thebibliography}
\end{document}